\documentclass[reqno]{amsart}
\usepackage[utf8]{inputenc}

\usepackage[utf8]{inputenc}

\usepackage{sky}
\usepackage{graphicx}
\usepackage{color}
\usepackage{caption}
\usepackage{subcaption}

\newcommand{\area}{\mrm{area}}

\newcommand{\vertiii}[1]{{\left\vert\kern-0.25ex\left\vert\kern-0.25ex\left\vert #1 
\right\vert\kern-0.25ex\right\vert\kern-0.25ex\right\vert}}
\newcommand{\sym}{S_{\mathrm{YM}}}

\title{Dynamical approach to area law for lattice Yang--Mills}

\author{Sky Cao}
\address{Department of Mathematics, Massachusetts Institute of Technology, Cambridge, MA 02139}
\email{skycao@mit.edu}

\author{Ron Nissim}
\email{rnissim@mit.edu}

\author{Scott Sheffield}
\email{sheffield@math.mit.edu}

\date{}

\begin{document}

\begin{abstract}
In this note, we prove Wilson's area law in the 't Hooft regime of parameters, which improves on a classical result of Osterwalder-Seiler from 1978, as well as on more recent work by the authors. The main point is to adapt the dynamical approach to lattice Yang--Mills set forth in \cite{shen2023stochastic} in order to verify the mass gap condition from \cite{durhuus1980connection}, from which area law directly follows. Our results apply for gauge groups $G \in \{\mathrm{U}(N), \mathrm{SU}(N), \mathrm{SO}(2N)\}$, which all have nontrivial center (which is one of the key assumptions in \cite{durhuus1980connection}). 
\end{abstract}

\maketitle

\tableofcontents

\section{Introduction}

The recent works \cite{shen2023stochastic, SZZ2024} laid out a stochastic analysis approach to the study of lattice Yang--Mills theories. Therein, it was observed that in the 't Hooft regime of parameters (to be introduced shortly), the Bakry-Emery condition may be verified, from which results such as exponential decay of correlations, the log-Sobolev and Poincar\'{e} inequalities, and uniqueness of the infinite volume limit all follow. However, one question that was not addressed in \cite{shen2023stochastic} was whether Wilson's area law would also follow from this approach. Wilson's area law is particularly important because it corresponds to {\em quark confinement} in the physics literature, see e.g.\ Wilson's original paper \cite{Wilson1974} or Chatterjee's recent survey \cite{chatterjee2016a}.

In a recent work \cite{CNS2025}, the authors managed via an entirely different approach to prove that area law holds in a certain regime of parameters (but not in the entire 't Hooft regime). This improved on a classical result of Osterwalder-Seiler \cite{osterwalderseiler1978}. The main point of the present paper is to obtain further improvements on \cite{osterwalderseiler1978, CNS2025}. Towards this end, we observe that in fact, area law in the 't Hooft regime follows quite readily from the dynamical approach. The additional input that one needs is the mass gap condition of \cite{durhuus1980connection}. We will discuss this condition in more detail in Section \ref{section:sigma-model-mass-gap-implies-area-law}. For now, we just roughly say that the idea is to split the lattice into a collection of height-$1$ slabs. Once boundary conditions on the top and bottom of a slab are fixed, what remains is to assign matrices to the edges that cross the slab, and the conditional law of these edges is given by a certain variable-environment $\sigma$-model. It is shown in \cite{durhuus1980connection} that if the resulting $\sigma$-models on each of these slabs exhibit exponential decay of correlations, uniform in the boundary conditions imposed at the top and bottom of the slabs, then area law follows.

\begin{remark}
There is also the recent work \cite{chatterjee2021probabilistic}, which shows that a different version of mass gap (see Definition 2.3 therein) implies area law. In principle, it may be possible to also verify this assumption starting from the Bakry-Emery condition, but we found it more direct to verify the condition of \cite{durhuus1980connection}.
\end{remark}

While \cite{shen2023stochastic} did not handle the case of $G = \mathrm{U}(N)$ (it focused on $G \in \{\mrm{SU}(N), \mrm{SO}(N)\}$), we observe that by a conditioning trick (essentially because $\mrm{U}(N) = \mrm{U}(1) \times \mrm{SU}(N)$ in a suitable sense), one may reduce the $\mrm{U}(N)$ case to the $\mrm{SU}(N)$ case, and thus our area law result applies also to $\mrm{U}(N)$. 

We now quickly introduce lattice Yang--Mills and then state the main result. Since we intend this to be a short note, we will be quite brief; see e.g.\ \cite{CPS2023, CNS2025, shen2023stochastic} for more discussion of the model. Let $\Lambda=\Lambda_L$ be the discrete torus of side length $L$, which can be written as $\Z^d$ modulo $L\Z^d$, and recall the Yang-Mills measure with gauge group $G \in \{\mathrm{SU}(N),\mathrm{U}(N),\mathrm{SO(N)}\}$ with 't Hooft scaling,
\begin{equs}\label{eq:Yang-Mills-measure}
    d\mu_{\mathrm{YM}}(Q) = \frac{1}{Z}\exp(\sym(Q))dQ,
\end{equs}
with $dQ = \prod_{e \in E_{\Lambda}^+} dQ_e $ the product Haar measure, and $\sym$ the Yang-Mills action defined by:
\begin{equ}\label{eq:sym}
    \sym(Q):= N\beta \sum_{p \in \mc{P}_{\Lambda}^+} \mathrm{Re} \mathrm{Tr}(Q_p).
\end{equ}
Here, $\mc{P}_\Lambda^+$ denotes the set of positively oriented plaquettes of $\Lambda$. Later on, we will also write $E_\Lambda^+$ for the set of positively oriented edges of $\Lambda$. To emphasize the dependence on $L$ and $\beta$, we sometimes refer to the measure as $\mu_{\mathrm{YM}_{L,\beta}}$.

\begin{remark}
In the above definitions and throughout this paper, $\Lambda = \Lambda_L$ is a torus, but this is not so essential.  For example, none of the arguments in this paper would substantially change if $\Lambda$ were instead a discrete cube of side length $L$ (with free boundary conditions). One should think of $L$ as arbitrarily large, and all our estimates will be uniform in $L$. 
\end{remark}

The main observables of interest for the model are the Wilson loop observables, which we define next. 

\begin{definition}[Loop variables and observables]
Let $Q : E_\Lambda^+ \ra G$ be a lattice gauge configuration. For a loop $\ell=e_1e_2\dots e_n$, let $Q_{\ell}:=Q_{e_1}Q_{e_2}\cdots Q_{e_n}$ denote the corresponding loop variable. We define the Wilson loop observable $W_\ell$ by
\begin{equs}
W_\ell(Q) := \tr(Q_\ell),
\end{equs}
where $\tr = \frac{1}{N} \Tr$ is the normalized trace. 
\end{definition}



Before stating the main results, we introduce a few $\beta$ thresholds which appear in the theorem statements.

\begin{definition}\label{def:beta-thresholds}
We define the following thresholds:
\begin{equs}
    \beta_{\mathrm{SU}(N)}^*=\beta_{\mathrm{U}(N)}^*&:= \frac{1}{8(d-1)},\\ \beta_{\mathrm{SO}(N)}^*&:= \frac{1}{16(d-1)}-\frac{1}{8N(d-1)}.
\end{equs}
\end{definition}

\begin{remark}
With the action $\sym$ as defined in \eqref{eq:sym} (in particular with the $N$ prefactor), the regime $\beta \leq \beta_0$ for some constant $\beta_0$ is known as the `t Hooft regime for lattice Yang--Mills theory. This is also the regime where one has a nontrivial large-$N$ limit, as studied in \cite{Chatterjee2019a, jafarov2016, chatterjee2016, BG2018, BCSK2024, BCSK25}, where $\beta_0$ is taken to be a small dimensional constant. 

The paper \cite{shen2023stochastic} also works in the 't Hooft regime, with the $\beta_0$ taken to be $ \beta^*_G / 2$ for $G \in \{\mrm{SU}(N), \mrm{SO}(N)\}$, i.e.\ the threshold appearing in \cite{shen2023stochastic} is half the threshold in the current paper. We obtain an improved threshold because we apply Bakry-Emery theory to a spin system defined on vertices, while \cite{shen2023stochastic} applies Bakry-Emery theory directly to lattice Yang--Mills, which is a spin system defined on edges. The combinatorial factors which appear in the proofs turn out to be a bit better for spin systems defined on vertices, essentially because every edge contains two vertices, while every plaquette contains four edges (and two is less than four). We thank Hao Shen, Rongchan Zhu, and Xiangchan Zhu for pointing this out to us.
\end{remark}

We now state the main result of the paper.

\begin{theorem}[Area law in the 't Hooft regime]\label{thm:main}
Let $d \geq 2$, $N \geq 2$, and $G \in \{\mrm{U}(N), \mrm{SU}(N), \mrm{SO}(2(N-1))\}$. Then for $\beta < \beta^*_G$, there are constants $C = C(\beta, d, N)$ and $c = c(\beta, d, N)$ such that for any rectangular loop $\ell$ in the lattice $\Lambda$, where the side lengths of $\ell$ are at most $L/2$, we have that
\begin{equ}
|\langle W_\ell \rangle_{\Lambda, \beta}| \leq C \exp(-c \,\, \area(\ell)).
\end{equ}
\end{theorem}

\begin{remark}
The restriction to $N \geq 2$ for $\mrm{SU}(N)$ arises because $\mrm{SU}(1) = \{1\}$ leads to a trivial model for which the result does not hold. The same restriction for $\mrm{U}(N)$ arises because our proof proceeds by reducing to the $\mrm{SU}(N)$ case, as previously mentioned. The restriction to even matrix sizes in $\mrm{SO}(2N)$ arises because in order to apply the result of \cite{durhuus1980connection}, the group $G$ must have a nontrivial center. This rules out the cases of $G = \mrm{SO}(2N+1)$. For more discussion on the role of the center and its relation to area law, see \cite{durhuus1980connection, chatterjee2021probabilistic}.

Since we intend for this paper to be a short note, we will not try to give an exhaustive list of references. For some previous work on Wilson's area law in various regimes, see \cite{Orland2005, Orland2006, Orland2007, Orland2008, Orland2020}. For a thorough review of the mathematical literature, see \cite[Section 6]{chatterjee2016a} 
\end{remark}

We now briefly summarize the rest of the paper. In Section \ref{section:sigma-model-mass-gap-implies-area-law}, we introduce the mass gap assumption of \cite{durhuus1980connection} which implies area law. In Section \ref{section:sigma-model-mass-gap}, we discuss how the Bakry-Emery condition implies the mass gap assumption of \cite{durhuus1980connection}, and additionally we use a conditioning trick to reduce the case $G = \mrm{U}(N)$ to the case $G = \mrm{SU}(N)$. This will give the proof of Theorem \ref{thm:main} for all cases of $G$. 

\begin{remark}
Some of our calculations are very similar to calculations given in \cite{shen2023stochastic}. In these cases, to avoid repeating the full arguments from \cite{shen2023stochastic}, we have indicated the pages of \cite{shen2023stochastic} where the calculations occur and have noted the few specific modifications needed for our setting. Specifically, we will invoke the proofs of Lemma 3.3, Lemma 4.1, and Corollary 4.11 from \cite{shen2023stochastic} (roughly seven pages of calculations). In order to follow the proofs in the current paper, the reader should have a copy of those arguments on hand. We also make heavy use of results from \cite{durhuus1980connection}, which we are able to cite more directly.
\end{remark}

\noindent \textbf{Acknowledgements:} We sincerely thank Martin Hairer, Peter Orland, Hao Shen, Rongchan Zhu, and Xiangchan Zhu for helpful discussions and comments. S.C. was partially supported by the NSF under Grant No. DMS-2303165. R.N. was supported by the NSF under Grant No. GRFP-2141064. S.S. was partially supported by the NSF under Grant No. DMS-2153742.

\section{\texorpdfstring{$\sigma$}{sigma}-Model Mass Gap Implies Area Law}\label{section:sigma-model-mass-gap-implies-area-law}

In this section, we discuss \cite[Theorems 1.2 and 1.3]{durhuus1980connection} using the notation of the present paper. The fundamental observation is that if we condition on all edges on the top and bottom horizontal boundaries of a given height-$1$ slab, then we obtain a $\sigma$-model with boundary conditions. In more detail, for some integer $k \in \Z$, we may condition on  the value of $Q_e$ for all edges $e$ lying in the ``horizontal'' planes $\{x=(x_1,\dots,x_d) \in \Z^d :x_d=k\}$ and $\{x=(x_1,\dots,x_d):x_d=k+1\}$. We write $Q_e = A_e$ for edges $e$ in the former plane, and $Q_e = B_e$ for edges in the latter plane. Having fixed the ``boundary conditions'' $A_e$ and $B_{e}$, we may further identify the ``vertical'' edges going between the two horizontal planes by their base point, and as a result the joint law of $(Q_e)_e$ for the vertical edges $e$ is given by the following nonlinear $\sigma$-model on vertices in $\Lambda^{d-1}=(\Z^d / L\Z^d)\times \{k\}$. 

\begin{definition}\label{def:Sigma-Model-intro}
Let $G \sse \mrm{U}(N)$ be a compact Lie group. Given boundary fields 
\begin{equ}
A=(A_e)_{e \in E_{\Lambda^{d-1}}^+},B=(B_e)_{e \in E_{\Lambda^{d-1}}^+} \in \mrm{U}(N)^{E_{\Lambda^{d-1}}^+},
\end{equ}
we define the action $S_{A,B}: G^{\Lambda^{d-1}} \to \R$ by 
\begin{equ}
S_{A,B}(Q):=N \beta\sum_{e=(x, y) \in E_{\Lambda^{d-1}}^+}\mathrm{Re}\mathrm{Tr}(Q_x A_e Q_y^{-1}B_e^{-1}),
\end{equ}
as well as the measure $\mu_{A, B}$ on $G^{\Lambda^{d-1}}$ by
    \begin{equs}\label{eq:sigma-model-intro}
        d\mu_{A,B}(Q) :=\frac{1}{Z_{A,B}} \exp(S_{A,B}(Q))dQ,
    \end{equs}
    where $dQ:=\prod_{v \in \Lambda^{d-1}} dQ_v$ is the product Haar measure, and $Z_{A,B}:= \int \exp(S_{A,B}(Q))dQ$.
\end{definition}

\begin{remark}
Here, we took the boundary fields $A, B$ to be $\mrm{U}(N)$-valued, even though $G \sse \mrm{U}(N)$ may not be all of $\mrm{U}(N)$. This will be convenient for later, when we apply the ``conditioning trick'' to extend our results from $\mrm{SU}(N)$ to $\mrm{U}(N)$. In particular, this is not needed if one wishes to just consider $\mrm{SU}(N), \mrm{SO}(2N)$.
\end{remark}

We now state a combined version of \cite[Theorems 1.2 and 1.3]{durhuus1980connection}.

\begin{theorem}[Theorems 1.2 and 1.3 of \cite{durhuus1980connection}]\label{thm:DF-mass-gap-implies-confinement}

Let $N \geq 1$ and $G \sse \mrm{U}(N)$ be a compact Lie group such that there exists $z \in \mrm{U}(1)$ with $z \neq 1, z I \in G$ (here $I$ is the $N \times N$ identity matrix). Let $d \geq 2$ and $\beta \geq 0$. For $i, j \in [N]$ and $x \in \Lambda^{d-1}$, let $f_{x}^{i, j}(Q)=(Q_x)_{ij}$ and $g_{x}^{i, j}(Q)=(Q_x^{-1})_{ij}$. Suppose that there exist constants $C_1$ and $C_2$ only depending on $G,d,\beta$ (and in particular not on $A, B$) such that for any $x,y \in \Lambda^{d-1}$ and $i_1, j_1, i_2, j_2 \in [N]$, we have that
\begin{equs}  \big|\mathrm{Cov}_{A,B}\big(f_x^{i_1 j_1}, g_y^{i_2 j_2}\big)\big| \leq C_1 e^{-C_2 d(x,y)},
\end{equs}
where the covariance is with respect to the $\sigma$-model $\mu_{A,B}$. Then the corresponding lattice Yang-Mills theory with the same $G$, $\beta$ and with dimension $d$ exhibits area law. I.e.\ there are constants $C_1, C_2$ only depending on $G, d, \beta$ such that for any rectangular loop $\ell$ in the lattice $\Lambda$, where the side lengths of $\ell$ are at most $L/2$, we have that
\begin{equs}\label{eq:Area-Law}
\big| \langle W_\ell \rangle_{\Lambda, \beta} \big| \leq C_1\exp(-C_2 \area(\ell)).
\end{equs}
%
\end{theorem}
\begin{proof}
We sketch the argument, referring the reader to \cite{durhuus1980connection} for a more complete explanation. Let $\ell_{R, T}$ be a rectangular loop with side lengths $R, T$. First, we condition on all the values of $Q_e$ orthogonal to the length-$T$ side of the rectangular loop. We refer to $T$ as the ``vertical side length''. Recalling that $W_{\ell_{R, T}} = \tr(Q_{\ell_{R, T}})$, we then expand $\mathrm{tr}(Q_{\ell_{R, T}})$ as a sum of $N^{|\ell_{R,T}|}$ terms which are products of $|\ell_{R,T}|$ matrix entries from matrices assigned to the edges along the loop. Since we conditioned all of the edges orthogonal to the vertical side of the rectangle, all the entries corresponding to ``horizontal'' edges on the length-$R$ side are fixed, and the remaining product of entries only depends on the height-$1$ slab to which they belong. Thus the terms in $\langle \mathrm{tr}(Q_{\ell_{R,T}})\rangle_{\Lambda, \beta}$ break down into products of $T$ two point correlation functions of the form 
\begin{equ}
\E_{A,B}[f_{x}^{i_1 j_1}g_{y}^{i_2 j_2}],
\end{equ}
where $d(x,y)=R$. Now by assumption, there exists $z \in \mrm{U}(1)$ such that $z \neq 1$, $z I \in G$. One may show that (refer to \cite[Theorem 1.3]{durhuus1980connection} for more details)
\begin{equ}
\E_{A, B}[f_x^{i_1 j_1}] = z \E_{A, B}[f_x^{i_1 j_1}].
\end{equ}
Thus since $z \neq 0$, the above expectation is zero, and thus
\begin{equ}
\E_{A,B}[f_{x}^{i_1 j_1}g_{y}^{i_2 j_2}] = \mrm{Cov}_{A, B}\big(f_{x}^{i_1 j_1}, g_{y}^{i_2 j_2}\big).
\end{equ}
By assumption, the right hand side is bounded (uniformly in $A, B$) by $C_1e^{-C_2 R}$. It thus follows that
\begin{equ}\label{eq:area-law-intermediate}
    \langle \mathrm{tr}(Q_{\ell_{R,T}}) \rangle_{\Lambda, \beta} \leq C_1N^{|\ell_{R,T}|}\exp(-C_2 \area(\ell_{R,T})).
\end{equ}
To prove \eqref{eq:Area-Law}, we recall that any lattice Yang--Mills theory with gauge group that has nontrivial center always (at least) exhibits the perimeter law (see e.g.\ \cite[Lemma 12.3]{chatterjee2021probabilistic}), i.e.\ 
\begin{equ}
\langle \mathrm{tr}(Q_{\ell_{R,T}})\rangle_{\Lambda, \beta} \leq C_1'\exp(-C_2' |\ell_{R,T}|).
\end{equ}
Thus if $\area(\ell_{R, T}) \geq C |\ell_{R, T}|$ for some appropriately chosen constant $C$, the perimeter prefactor $N^{|\ell_{R, T}|}$ in \eqref{eq:area-law-intermediate} may be absorbed into the area decay and we obtain area law, otherwise if $\area(\ell_{R, T}) < C|\ell_{R, T}|$, then the perimeter and area are comparable, and so the above perimeter law is enough to imply an area law bound.
\end{proof}

\section{\texorpdfstring{$\sigma$}{sigma}-Model Mass Gap}\label{section:sigma-model-mass-gap}

In this section, we observe that the arguments of \cite{shen2023stochastic} can be used to verify the assumptions of Theorem \ref{thm:DF-mass-gap-implies-confinement}. We first discuss the cases $G \in \{\mrm{SU}(N), \mrm{SO}(N)\}$, before moving on to the case $G = \mrm{U}(N)$. Recall that $\Lambda^{d-1}$ denotes a $(d-1)$-dimensional slice of the $d$-dimensional lattice $\Lambda=\Lambda_L$.

\begin{notation}[Local observables]
A local observable is a function $f \colon G^{\Lambda_f} \ra \R$, where $\Lambda_f \sse \Z^{d-1}$ is a finite subset of vertices. In particular, whenever we specify a local observable $f$, we implicitly also specify the set $\Lambda_f$. Additionally, for any lattice $\tilde{\Lambda} \sse \Z^{d-1}$ such that $\Lambda_f \sse \tilde{\Lambda}$, we have may extend $f$ to a function $f \colon G^{\tilde{\Lambda}} \ra \R$. Given a local observable $f$, let
\begin{equs}
    \vertiii{f}_{\infty} := \sum_{x \in \Lambda_f} \| \nabla_x f\|_{L^{\infty}}.
\end{equs}
\end{notation}

We now state the first main result of the section, which is that the assumptions of Theorem \ref{thm:DF-mass-gap-implies-confinement} hold for $G \in \{\mrm{SU}(N), \mrm{SO}(N)\}$.

\begin{prop}[Uniform Mass Gap for $\sigma$-Model]\label{prop:uniform-mass-gap-SU(N)}
Let $N \geq 2$. For $G \in \{\mathrm{SU}(N),\mathrm{SO}(N)\}$, $\beta < \beta_G^*$, and any choice of fields $A,B \in \mrm{U}(N)^{E_{\Lambda^{d-1}}^+}$, there are constants $C_1,C_2$ only depending on $G, d, \beta$, such that for any local observables $f,g \in C^{\infty}(G^{\Lambda^{d-1}})$, we have that
    \begin{equs}
        \big| \mathrm{Cov}_{A, B}(f,g)\big| \leq C_1 e^{-C_2 d(\Lambda_f, \Lambda_g)}(\vertiii{f}_{\infty} \vertiii{g}_{\infty} + \|f\|_{L^2(\mu_{A, B})}\|g\|_{L^2(\mu_{A, B})}).
    \end{equs}
    Here, the covariance is with respect to $\mu_{A, B}$ (recall Definition \ref{def:Sigma-Model-intro}).
\end{prop}


Before proving Proposition \ref{prop:uniform-mass-gap-SU(N)}, we make the trivial observation that Theorem \ref{thm:DF-mass-gap-implies-confinement} now follows for $G \in \{\mrm{SU}(N), \mrm{SO}(N)\}$.

\begin{proof}[Proof of Theorem \ref{thm:main} for $G \in \{\mrm{SU}(N), \mrm{SO}(N)\}$.]
This now follows by combining Theorem \ref{thm:DF-mass-gap-implies-confinement} with Proposition \ref{prop:uniform-mass-gap-SU(N)}.
\end{proof}

\begin{proof}[Proof of Proposition \ref{prop:uniform-mass-gap-SU(N)}]
For an introduction to the geometric preliminaries used in the proof, see \cite[Section 2]{shen2023stochastic}. Similarly, as much of the argument is the exact same as in \cite{shen2023stochastic}, we leave many of the details to the cited reference, although we will provide precise citations to the corresponding locations in \cite{shen2023stochastic}. First, by exactly the same arguments as in the proof of \cite[Lemma 4.1]{shen2023stochastic}, we may show that for any boundary fields $A, B \in G^{E_{\Lambda^{d-1}}^+}$, any field $Q \in G^{\Lambda^{d-1}}$, and any tangent vector $v = X Q \in T_Q G^{\Lambda^{d-1}}$, we have the Hessian bound (recall Definition \ref{def:Sigma-Model-intro} for the definition of $S_{A, B}$)
\begin{equ}\label{eq:hessian-bd}
|\mrm{Hess}_{S_{A, B}}(v, v)| \leq 4(d-1) N \beta |v|^2.
\end{equ}
In more detail, referring to the top of \cite[Page 823]{shen2023stochastic} for the start of the proof of \cite[Lemma 4.1]{shen2023stochastic}, one sees that even with the additional boundary conditions $A, B$, one may perform the exact same computations as in \cite[(4.3)-(4.4)]{shen2023stochastic} and the next three math displays thereafter. Then in the fourth math display after \cite[(4.4)]{shen2023stochastic}, one applies Cauchy-Schwarz, and the additional boundary conditions $A, B$ disappear completely (due to the fact that $A, B$ consist of unitary matrices). To be precise, one should replace $E_\Lambda^+$ therein by our $\Lambda^{d-1}$. From this, one obtains the intermediate estimate
\begin{equ}
\big|\mrm{Hess}_{S_{A, B}}(v, v)\big| \leq N \beta \sum_{x \in \Lambda^{d-1}} \sum_{\substack{e \in E_{\Lambda^{d-1}}^+ \\ e \ni x}}
|X_x|^2 + N \beta \sum_{\substack{x \neq y \in \Lambda^{d-1} \\ x \sim y}} \frac{1}{2} \big(|X_x|^2 + |X_{x'}|^2\big),
\end{equ}
where here $e \ni x$ denotes that the edge $e$ contains the vertex $x$, and $x \sim y$ denotes that $(x, y)$ is an (oriented) edge of $\Lambda^{d-1}$. The first term on the right hand side corresponds to the term estimated in fifth math display after \cite[(4.4)]{shen2023stochastic}, while the second term corresponds to the term estimated in the sixth math display after \cite[(4.4)]{shen2023stochastic}. To further estimate the above, we use that each vertex of $\Lambda^{d-1}$ is contained in $2(d-1)$ edges of $\Lambda^{d-1}$, from which we obtain
\begin{equs}
\big|\mrm{Hess}_{S_{A, B}}(v, v)\big| &\leq 2(d-1) N \beta \sum_{x \in \Lambda^{d-1}} |X_x|^2 + N \beta \sum_{x \in \Lambda^{d-1}} |X_x|^2 \sum_{e \ni x} 1 \\
&\leq 4(d-1) N \beta |v|^2,
\end{equs}
where here $|v|$ is the norm of the tangent vector $v$, which has the more explicit form
\begin{equ}
|v|^2 = \sum_{x \in \Lambda^{d-1}} |X_x|^2.
\end{equ}
This shows the Hessian bound \eqref{eq:hessian-bd}. 
From this, it follows that 
the Bakry-Emery condition holds:
\begin{equ}\label{eq:Bakry-Emery}
\mrm{Ric}_{G^{\Lambda^{d-1}}} (v, v) - \mrm{Hess}_{S_{A, B}}(v, v) \geq K_{S_{A, B}} |v|^2,
\end{equ}
where $K_{S_{A, B}}$ is almost exactly the same as $K_{\mc{S}}$ in \cite[Assumption 1.1]{shen2023stochastic}, i.e.
\begin{equ}
K_{S_{A, B}} = \begin{cases} \frac{N+2}{4} - 1 - 4N\beta (d-1), & G = \mrm{SO}(N), \\
\frac{N+2}{2} - 1 - 4N\beta(d-1) , & G = \mrm{SU}(N).\end{cases}
\end{equ}
The argument establishing \eqref{eq:Bakry-Emery} is almost exactly the same as in the proof of \cite[Theorem 4.2]{shen2023stochastic}. The only difference is that our $4N \beta(d-1)$ is $8N\beta(d-1)$ in \cite[Assumption 1.1]{shen2023stochastic}, which is due to the fact that the Hessian bound proven in \cite[Lemma 4.1]{shen2023stochastic} involves a prefactor 8 (whereas we have a prefactor 4 in \ref{eq:hessian-bd}). Similar to the remark right after \cite[Assumption 1.1]{shen2023stochastic}, the assumption that $\beta < \beta_G^*$ (recall Definition \ref{def:beta-thresholds}) is equivalent to $K_{S_{A, B}} > 0$.
Then as in \cite[Corollary 4.4 and Remark 4.6]{shen2023stochastic}, the Bakry-Emery condition \eqref{eq:Bakry-Emery} implies that for any local observable $f$,
\begin{equ}\label{eq:be-variance-decay}
\big(\mrm{Var}_{A, B}(P_t f)\big)^{1/2} \leq e^{- K_{S_{A, B}} t} \|f\|_{L^2(\mu_{A, B})},
\end{equ}
where $P_t$ is a semigroup with generator $\mc{L}$ given by
\begin{equ}\label{eq:Pt-generator}
\mc{L} F(Q) = \sum_{x \in V_\Lambda} \Delta_x F(Q) + \sum_{x \in V_\Lambda} \langle \nabla_x S_{A, B}(Q), \nabla_x F \rangle.
\end{equ}
The remainder of the argument exactly follows the proof of \cite[Corollary 4.11]{shen2023stochastic}, thus we only sketch the ideas here. First, by the invariance of $\mu_{A, B}$ with respect to the semigroup $P_t$ (see e.g.\ \cite[Lemma 3.3]{shen2023stochastic}), we may write
\begin{equ}\label{eq:covariance-decomposition-intermediate}
\mrm{Cov}_{A, B} (f, g) = \E_{A, B} \big(P_t(fg) - P_t f P_t g\big) + \mrm{Cov}_{A, B}(P_t f, P_t g).
\end{equ}
Take $t \sim d(\Lambda_f, \Lambda_g)$. The latter term may be directly estimated by \eqref{eq:be-variance-decay} (see e.g.\ \cite[(4.27)]{shen2023stochastic}). To estimate the former term, first note that when $t = 0$, $P_t (fg) - P_t f P_t g = 0$. Next, since $f, g$ are local observables with supports separated at distance $d(\Lambda_f, \Lambda_g)$, and since $P_t$ is a semigroup whose generator has a ``local" structure (which is reflected in the fact that the action $S_{A, B}$ is defined by only looking at neighboring spins), one expects that for $t$ not too large, we still have that $P_t (fg) - P_t f P_t g$ is approximately zero. More quantitatively, we expect that if $t \sim c d(\Lambda^{d-1}_f, \Lambda^{d-1}_g)$ for some sufficiently small constant $c$, then 
\begin{equ}
|P_t(fg) - P_t f P_t g| \leq C e^{-c d(\Lambda_f, \Lambda_g)} \vertiii{f}_{L^\infty} \vertiii{g}_{L^\infty}.
\end{equ}
This is indeed shown in the proof of \cite[Corollary 4.11]{shen2023stochastic}\footnote{From \cite[Remark 4.12]{shen2023stochastic}, the constant $c$ that one obtains decays with $N$, but this is not necessarily optimal, as mentioned in the cited remark.}; see the argument appearing after \cite[(4.27)]{shen2023stochastic}, which remains completely unchanged for us. Combining the estimates for the two terms in \eqref{eq:covariance-decomposition-intermediate}, the covariance bound now follows.
\end{proof}

Having proven Theorem \ref{thm:main} in the cases $G \in \{\mrm{SU}(N), \mrm{SO}(N)\}$, we now turn to the case $G = \mrm{U}(N)$. Because $\mathrm{U}(N)$ does not have uniformly positive Ricci curvature, one cannot directly verify the Bakry-Emery condition in this case. Instead, we proceed by using the fact that $\mrm{U}(N) = \mrm{U}(1) \times \mrm{SU}(N)$, in the sense that the product of Haar-distributed $\mrm{U}(1)$ and $\mrm{SU}(N)$ matrices is a Haar-distributed $\mrm{U}(N)$ matrix. Thus, we may condition on the $\mrm{U}(1)$ part, leaving the $\mrm{SU}(N)$ part, which we already know how to handle. Besides this conditioning trick, we also need to introduce a few more definitions before stating the main proposition for the $\mrm{U}(N)$ case.


\begin{definition}\label{def:U-N-SU-N-function}
    For any function $f: \mathrm{U}(N)^{\Lambda^{d-1}} \to \C$ and some field $(z_x)_{x \in \Lambda^{d-1}} \in \mathrm{U}(1)^{\Lambda^{d-1}}$, we define $f_z(\tilde{Q})=f(z\tilde{Q})$ as a function on $\mathrm{SU}(N)^{\Lambda^{d-1}}$.
\end{definition}

\begin{definition}
We say that a function $f : \mathrm{U}(N)^{\Lambda^{d-1}} \ra \C$ is linear if 
\begin{equ}
f(e^{i \theta} Q) = e^{i \theta} f(Q), ~~ Q \in \mrm{U}(N)^{\Lambda^{d-1}}, ~~ \theta \in \R,
\end{equ}
where $e^{i \theta} Q := (e^{i \theta}Q_x)_{x \in \Lambda^{d-1}}$. This is a slight abuse of terminology, because $\mrm{U}(N)^{\Lambda^{d-1}}$ is not a vector space.
\end{definition}

\begin{remark}\label{remark:linear-assumption}
We remark that the function function $f_x^{i_1 j_1}$ appearing in Theorem \ref{thm:DF-mass-gap-implies-confinement} is linear, as it is defined as the $(i_1, j_1)$ matrix entry of the matrix $Q_x$ at the vertex $x$.
\end{remark}

We now state and prove the following corollary of Proposition \ref{prop:uniform-mass-gap-SU(N)}.

\begin{cor}[Uniform mass gap for $\sigma$-Model, $\mathrm{U(N)}$ Case]\label{cor:uniform-mass-gap-U(N)}
Let $N \geq 2$ and $G = \mathrm{U}(N)$. There are constants $C_1,C_2$ only depending on $G, d, \beta$ such that the following holds for $\beta < \beta_{\mrm{U}(N)}^*$. For any choice of fields $A, B \in G^{E_{\Lambda^{d-1}}^+}$
and any local observables $f,g \in C^{\infty}(G^{\Lambda^{d-1}})$ such that $f$ is linear, we have that
\begin{equs}
\big| \mathrm{Cov}_{A, B}(f,g) \big| \leq C_1 e^{-C_2 d(\Lambda_f, \Lambda_g)}\sup_{z \in \mathrm{U}(1)^{\Lambda^{d-1}}}\big(\vertiii{f_z}_{\infty} \vertiii{g_z}_{\infty} + \|f_z\|_{L^\infty}\|g_z\|_{L^\infty}\big),
\end{equs}
where the covariance is with respect to $\mu_{A, B}$.
\end{cor}
\begin{proof}
Taking two independent random variables $z,\tilde{Q}$ such that  $\tilde{Q} \sim \mathrm{Haar}(\mathrm{SU}(N))$, and $z \sim \mathrm{Haar}(\mathrm{U}(1))$, we have that $z \tilde{Q} \sim \mathrm{Haar}(\mathrm{U}(N))$. This allows us to split Haar integration over $\mrm{U}(N)$ into a product of Haar integration over $\mrm{U}(1)$ and Haar integration over $\mathrm{SU}(N)$. In particular, for any observable $f$, we have that
    \begin{equs}
        \int\prod_{v \in \Lambda^{d-1}}dQ_v\, f(Q) \exp\big(S_{A, B}(Q)\big) &= \int \prod_{v \in \Lambda^{d-1}} dz_v \int \prod_{v \in \Lambda^{d-1}}d\tilde{Q}_v \,  f(z\tilde{Q})  \exp\big(S_{\tilde{A}, B}(\tilde{Q})\big) , \\
        &= \int \prod_{v \in \Lambda^{d-1}} dz_v Z_{\tilde{A}, B} \int  f_z(\tilde{Q}) d\mu_{\tilde{A}, B}(\tilde{Q}), \label{eq:disintegration-identity}
    \end{equs}
    where $\tilde{A}_e = z_x z_y^{-1}A_e$, and $S_{\tilde{A}, B}, \mu_{\tilde{A}, B}, Z_{\tilde{A}, B}$ are as in Definition \ref{def:Sigma-Model-intro} (with the $G$ therein taken to be $\mrm{SU}(N)$). Applying this identity with $f \equiv 1$, we also have the normalizing constant identity
    \begin{equ}\label{eq:partition-fn-integration-identity}
    Z_{A, B} = \int \prod_{v \in \Lambda^{d-1}} dQ_v\,  \exp(S_{A, B}(Q))  = \int \prod_{v \in \Lambda^{d-1}} dz_v \, Z_{\tilde{A}, B}.
    \end{equ}
    We will shortly use the following result: for fixed $z = (z_v)_{v \in \Lambda^{d-1}}$, applying the $\mathrm{SU}(N)$ uniform mass gap result (Proposition \ref{prop:uniform-mass-gap-SU(N)}) (and further bounding the $L^2(\mu_{\tilde{A}, B})$ norm by the $L^\infty$ norm), we have that
    \begin{equ}\label{eq:su(n)-uniform-mass-gap-application}
    \big|\Cov_{\tilde{A}, B}(f_z, g_z)\big| \leq  C_1 e^{-C_2 d(\Lambda_f, \Lambda_g)}\big(\vertiii{f_z}_{\infty} \vertiii{g_z}_{\infty} + \|f_z\|_{L^\infty}\|g_z\|_{L^\infty}\big).
    \end{equ}
    Next, since $N \geq 2$, we may take $k \in [N-1]$ such that $e^{i 2\pi k / N} \neq 1$. Given $\tilde{Q} \in \mrm{SU}(N)^{\Lambda^{d-1}}$, let $\tilde{Q}' = e^{i 2\pi k / N} Q$. The action $S_{A, B}$ (recall Definition \ref{def:Sigma-Model-intro}) remains invariant under this multiplication, i.e.\ $S_{\tilde{A}, B}(\tilde{Q}) = S_{\tilde{A}, B}(\tilde{Q}')$. Additionally, Haar measure on $\mrm{SU}(N)$ remains invariant under this multiplication. Then since $f$ is linear by assumption, it follows that
    \begin{equs}
    \int f_z(\tilde{Q}) d\mu_{\tilde{A}, B}(\tilde{Q}) = \int f_z(\tilde{Q}') d\mu_{\tilde{A}, B}(\tilde{Q}) = e^{i 2\pi k / N} \int f_z(\tilde{Q}) d\mu_{\tilde{A}, B}(\tilde{Q}),
    \end{equs}
    and since $e^{i 2\pi k / N} \neq 1$, it follows that the above is zero. Combining this with \eqref{eq:disintegration-identity}, \eqref{eq:su(n)-uniform-mass-gap-application}, and \eqref{eq:partition-fn-integration-identity}, we obtain
    \begin{equs}
    \big|\mathrm{Cov}_{A, B} (f,g)\big| &= \bigg| Z_{A, B}^{-1} \int \prod_{v \in \Lambda^{d-1}} dz_v \, Z_{\tilde{A}, B} \Cov_{\tilde{A}, B}(f_z, g_z)\bigg| \\
    &\leq \bigg(Z_{A, B}^{-1} \int \prod_{v \in \Lambda^{d-1}} dz_v \, Z_{\tilde{A}, B} \bigg) C_1 e^{-C_2 d(\Lambda_f, \Lambda_g)}\sup_{z \in \mathrm{U}(1)^{\Lambda^{d-1}}}\big(\vertiii{f_z}_{\infty} \vertiii{g_z}_{\infty} + \|f_z\|_{L^\infty}\|g_z\|_{L^\infty} \big), \\
    &= C_1 e^{-C_2 d(\Lambda_f, \Lambda_g)}\sup_{z \in \mathrm{U}(1)^{\Lambda^{d-1}}} \big(\vertiii{f_z}_{\infty} \vertiii{g_z}_{\infty} + \|f_z\|_{L^\infty}\|g_z\|_{L^\infty} \big),
    \end{equs} 
as desired.
\end{proof}

\begin{proof}[Proof of Theorem \ref{thm:main} for $G = \mrm{U}(N)$.]
This now follows by combining Theorem \ref{thm:DF-mass-gap-implies-confinement} and Corollary \ref{cor:uniform-mass-gap-U(N)}, where we use the fact that $f_x^{i_1 j_1}$ is linear (recall Remark \ref{remark:linear-assumption}).
\end{proof}

\bibliographystyle{alpha}
\bibliography{references}

\begin{thebibliography}{Cha19b}

\bibitem[BCS24]{BCSK2024}
Jacopo {Borga}, Sky {Cao}, and Jasper {Shogren-Knaak}.
\newblock {Surface sums for lattice Yang-Mills in the large-$N$ limit}.
\newblock {\em arXiv e-prints}, page arXiv:2411.11676, November 2024.

\bibitem[BCS25]{BCSK25}
Jacopo {Borga}, Sky {Cao}, and Jasper {Shogren-Knaak}.
\newblock {Surface sums in two-dimensional large-$N$ lattice Yang--Mills: Cancellations and explicit computation for general loops}.
\newblock {\em arXiv e-prints}, page arXiv:2508.13827, August 2025.

\bibitem[BG18]{BG2018}
Riddhipratim Basu and Shirshendu Ganguly.
\newblock {$\mathrm{S}\mathrm{O}(N)$} lattice gauge theory, planar and beyond.
\newblock {\em Commun. Pure Appl. Math.}, 71(10):2016--2064, 2018.

\bibitem[Cha19a]{Chatterjee2019a}
Sourav Chatterjee.
\newblock Rigorous solution of strongly coupled {SO($N$)} lattice gauge theory in the large {$N$} limit.
\newblock {\em Communications in Mathematical Physics}, 366(1):203--268, 2019.

\bibitem[Cha19b]{chatterjee2016a}
Sourav Chatterjee.
\newblock Yang-{Mills} for probabilists.
\newblock In {\em Probability and Analysis in Interacting Physical Systems}, pages 1--16. Springer, Cham: Springer, 2019.

\bibitem[Cha21]{chatterjee2021probabilistic}
Sourav Chatterjee.
\newblock A probabilistic mechanism for quark confinement.
\newblock {\em Communications in Mathematical Physics}, 385(2):1007--1039, 2021.

\bibitem[CJ16]{chatterjee2016}
Sourav Chatterjee and Jafar Jafarov.
\newblock The {$1/N$} expansion for {SO($N$)} lattice gauge theory at strong coupling.
\newblock {\em arXiv preprint arXiv:1604.04777}, 2016.

\bibitem[CNS25]{CNS2025}
Sky {Cao}, Ron {Nissim}, and Scott {Sheffield}.
\newblock {Expanded regimes of area law for lattice Yang-Mills theories}.
\newblock {\em arXiv e-prints}, page arXiv:2505.16585, May 2025.

\bibitem[CPS23]{CPS2023}
Sky {Cao}, Minjae {Park}, and Scott {Sheffield}.
\newblock {Random surfaces and lattice Yang-Mills}.
\newblock {\em arXiv e-prints}, page arXiv:2307.06790, July 2023.
\newblock To appear in {\em Comm. Amer. Math. Soc.}

\bibitem[DF80]{durhuus1980connection}
B~Durhuus and J~Fr{\"o}hlich.
\newblock A connection between $\nu$-dimensional yang-mills theory and ($\nu$- 1)-dimensional, non-linear $\sigma$-models.
\newblock {\em Communications in Mathematical Physics}, 75(2):103--151, 1980.

\bibitem[Jaf16]{jafarov2016}
Jafar Jafarov.
\newblock Wilson loop expectations in {SU($N$)} lattice gauge theory.
\newblock {\em arXiv preprint arXiv:1610.03821}, 2016.

\bibitem[Orl05]{Orland2005}
Peter Orland.
\newblock $(2+1)$-dimensional lattice qcd.
\newblock {\em Phys. Rev. D}, 71:054503, Mar 2005.

\bibitem[Orl06]{Orland2006}
Peter Orland.
\newblock Integrable models and confinement in $(2+1)$-dimensional weakly-coupled yang-mills theory.
\newblock {\em Phys. Rev. D}, 74:085001, Oct 2006.

\bibitem[Orl07]{Orland2007}
Peter Orland.
\newblock String tensions and representations in anisotropic ($2+1$)-dimensional weakly-coupled yang-mills theory.
\newblock {\em Phys. Rev. D}, 75:025001, Jan 2007.

\bibitem[Orl08]{Orland2008}
Peter Orland.
\newblock Composite strings in ($2+1$)-dimensional anisotropic weakly coupled yang-mills theory.
\newblock {\em Phys. Rev. D}, 77:025035, Jan 2008.

\bibitem[Orl20]{Orland2020}
Peter Orland.
\newblock {Confinement for all couplings in a ${\mathbb Z}_{2}$ lattice gauge theory}.
\newblock {\em J. Phys. A}, 53(13):13LT01, 2020.

\bibitem[OS78]{osterwalderseiler1978}
K~Osterwalder and E~Seiler.
\newblock Gauge field theories on a lattice.
\newblock {\em Annals of Physics}, 110(2):440--471, 1978.

\bibitem[SZZ23]{shen2023stochastic}
Hao Shen, Rongchan Zhu, and Xiangchan Zhu.
\newblock A stochastic analysis approach to lattice {Y}ang-{M}ills at strong coupling.
\newblock {\em Comm. Math. Phys.}, 400(2):805--851, 2023.

\bibitem[SZZ24]{SZZ2024}
Hao {Shen}, Rongchan {Zhu}, and Xiangchan {Zhu}.
\newblock {Langevin dynamics of lattice Yang-Mills-Higgs and applications}.
\newblock {\em arXiv e-prints}, page arXiv:2401.13299, January 2024.

\bibitem[Wil74]{Wilson1974}
Kenneth~G Wilson.
\newblock Confinement of quarks.
\newblock {\em Physical review D}, 10(8):2445, 1974.

\end{thebibliography}

\end{document}